 \documentclass[letterpaper, onecolumn]{IEEEtran}

\usepackage{graphicx,amstext,amsmath,amssymb,psfrag,enumerate,color}
\usepackage{parskip,upgreek, times}
\usepackage{soul,here}

\newtheorem{theorem}{Theorem}
\newtheorem{lemma}[theorem]{Lemma}
\newtheorem{definition}[theorem]{Definition}
\newtheorem{proposition}[theorem]{Proposition}

\newtheorem{assumption}{Assumption}
\newtheorem{remark}{Remark}

\graphicspath{{figures_recentes/}}

\def\0{{\bf 0}}
\def\R{\mathbb{R}}
\def\N{\mathbb{N}}

\def\Q{\mathbb{Q}}

 \def\V{\mathcal{V}}
  \def\S{\mathcal{S}}

\def\G{\mathcal{G}}

\def\E{\mathcal{E}}

\newcommand{\ie}{{\it i.e. }}

\usepackage{xcolor}
\usepackage{lipsum}
\newcommand{\blue}[1]{\textcolor{black}{#1}}
\newcommand{\sam}[1]{\blue{#1}}

\begin{document}

\title{ Continuous opinions and discrete actions in social networks: a multi-agent system approach\thanks{The work of I.-C. Mor\u{a}rescu, S. Martin  was funded by the ANR project COMPACS - "Computation Aware Control Systems", ANR-13-BS03-004.} }

\author{N. R. Chowdhury, I.-C. Mor\u{a}rescu, S. Martin, S. Srikant \thanks{I.-C. Mor\u{a}rescu and S. Martin  are with Universit\'e de Lorraine, CRAN, UMR 7039 and CNRS, CRAN, UMR 7039, 2 Avenue de la For\^et de Haye, Vand\oe uvre-l\`es-Nancy, France.  e-mails: {\tt  constantin.morarescu@univ-lorraine.fr, samuel.martin@univ-lorraine.fr}}}

\maketitle

\begin{abstract}
This paper proposes and analyzes a novel multi-agent opinion dynamics model in which agents have access to actions which are quantized version of the opinions of their neighbors. The model produces different behaviors observed in social networks such as disensus, clustering, oscillations, opinion propagation, \sam{even when the communication network is connected}. The main results of the paper provides the characterization of preservation and diffusion of actions under general communication topologies. A complete analysis allowing the opinion forecasting is given in the particular cases of complete and ring communication graphs. Numerical examples illustrate the main features of this model.
\end{abstract}

\IEEEpeerreviewmaketitle

\section{Introduction}\label{sec:intro}

The analysis of social networks received an increasing interest during the last decade. This is certainly related with the increasing use of Facebook, Twitter, LinkedIn and other on-line platforms allowing to get information about social networks. A \sam{convenient} way to model opinion dynamics in social networks is by the intermediate of multi-agent systems. The existing models can be split in two main classes: those considering that opinions can evolve in a discrete set and those considering a continuous set of values that can be taken by each agent. The models in first class come from statistical physics and the most employed are the Ising~\cite{IsingThesis1924}, voter \cite{CliffordSudbury1973} and Sznajd \cite{Sznajd2000} models.  When the opinions are not constrained to \sam{a discrete set}, we can find in the literature two popular models: the Deffuant \cite{Deffuant2000} and the Hegselmann-Krause~\cite{krause2002} models. They are usually known as bounded confidence models since they depend on one parameter characterizing the fact that one agent takes into account the opinion of another only if their opinions are close enough.
\sam{The bounded confidence models above do not guaranty consensus and instead several local agreements can be reached. The  Hegselmann-Krause model has been adapted in \cite{MG10} to a model of opinion dynamics with decaying confidence.}

In order to more accurately describe the opinion dynamics and to recover more realistic behaviors, a mix of continuous opinion with discrete actions (CODA) was proposed in \cite{Martins2008}. This model reflects the fact that even if we often face binary \sam{choices or} actions, the opinion behind evolves in a continuous space of values. For instance we may think that car A is 70\% more appropriate for our use than car B. However, our action will be 100\% buy car A. Moreover, our neighbors \sam{often} only see our action without any access to our opinion. Similar idea was employed in \cite{CeragioliFrascaECC2015} where it is studied the emergence of consensus under quantized all-to-all communication. In this paper the authors assume constant interaction weights and quantized information on the opinion of all the other individuals belonging to a given social network. \sam{In~\cite{frasca2012continuous}, a system with quantized communication and general communication topologies is studied. There, the author focuses on consensus (up to the quantizer precision) while in the present paper, we precisely focus on dynamics such as dissensus which occur within the quantizer precision. As explained abobve, this is of particular relevance for the social sciences.}

Whatever is the model employed to describe the opinion dynamics, many studies focus on the emergence of consensus in social networks \cite{GalamMoscovici1991,Axelrod1997,Deffuant2000, Fortunato2005}. Nevertheless, this behavior is rarely observed in real large social networks. \sam{The present study provides a possible explanation for dissensus even though the communication graph is connected.}

From mathematical view-point the model proposed in this paper is close to the one in \cite{CeragioliFrascaECC2015} with the difference that we are considering state-dependent interaction weights instead of constant ones. Beside the heterogeneity introduced in the model by the state-dependent interaction weights we are also dealing with more general interaction topologies and we are not trying to characterize consensus. Instead we highlight that extremist individuals are less influenceable, that several equilibria can be reached and one can also have oscillatory behaviors in the network.  From behavioral view-point our model is close to the one in \cite{Martins2008} \sam{(this similarity is highlighted in Section~\ref{sec:martins})}. The difference here is that we are able to analytically study this model and go beyond simulations and theirs interpretation.

The main contributions of this work are as follows:
\begin{itemize}
\item {\bf we propose a consensus-like dynamics} that approaches the dynamics described in   \cite{Martins2008}. This CODA model is given by a quantized consensus system with state-dependent interaction weights.

\item {\bf we describe the possible equilibria} of the proposed model and {\bf depict the main properties} characterizing the opinions dynamics. Precisely, we provide a criterion to detect stabilization of the actions of a group of agents and to predict the propagation of this action throughout the network. Our criterion depends on the initial conditions and interaction topology only.

\item {\bf we completely analyze} some particular interaction topologies such as: the complete graph and the ring graph.
\end{itemize}

The rest of the paper is organized as follows. Section \ref{problem_formulation} formulates the problem and illustrates that our model is close to the one proposed in \cite{Martins2008}. In Section \ref{COCAmodel} we show when the quantization effect is removed \ie in the context of continuous opinions with continuous actions (COCA), that consensus is always achieved provided that the interaction graph is connected. The main features of our CODA model under general interaction topologies are derived in Section \ref{CODAanalysis}. Precisely we characterize the preservation and the propagation of actions inside a group and outside it, respectively. Some particular interaction topologies such as: the complete graph and the ring graph are analyzed in \ref{Particular-Top}. 

\section{Problem formulation and preliminaries}\label{problem_formulation}

\subsection{Model description}

In this work we consider a network of $n$ individuals/agents denoted by $\V=\{1,\ldots, n\}$. The interaction topology between agents is described by a fixed graph $\G=(\V,\E)$ that can be directed or not. Let us also denote by $N_i$ the set of agents that influence $i$ according to the graph $\G$ (\ie $j\in N_i \Leftrightarrow (j,i)\in\E$) and by $n_i$ the cardinality of $N_i$. Initially, agent $ i\in\V $ has a given opinion $ p_i(0)=p_i^0 \in[0,1]$ and this opinion evolves according to a discrete-time protocol. Let $p_i(k)$ be the opinion of agent $i\in\V$ at time $k$. Assume that $\forall i\in\V, \ p_i(0)\neq\frac{1}{2}$ and introduce the action value $q_i(k)\in\{0,1\}$ as a quantized version of $p_i(k)$ defined by:

\begin{equation*}
q_i(k)= \left\{
\begin{array}{l}
0  \text{ if }\left( p_i(k)<\frac{1}{2}\right),\\
0  \text{ if }\left( p_i(k)=\frac{1}{2}~\text{and}~p_i(k-1)<\frac{1}{2}\right),\\
 1  \text{ otherwise.}
\end{array}\right.
\end{equation*}

Two distinct situations are considered in the following.\\
{\bf COCA model:} each agent has access to the opinion of its neighbors. In this case, our model simply writes as a consensus protocol with state-dependent interaction weights. Precisely, the opinion of an agent $i\in\V$ updates according to the following rule
\begin{equation}
p_i(k+1)=p_i(k)+\frac{p_i(k)(1-p_i(k))}{n_i}\sum_{j\in N_i} (p_j(k)-p_i(k)). \label{eq:mn1}
\end{equation}

Denoting by $p(\cdot)=(p_1(\cdot),\ldots,p_n(\cdot))$ the vector that collects the opinions of all agents, the collective opinion dynamics is given by:
\begin{equation}
p(k+1)=\left( I_n+A(p)\right) p(k).   \label{eq:mn2_vector}
\end{equation}

\begin{remark}
We assume that $p_i^0$ belongs to $(0,1)$ which is a normalized version of $\R$. Doing so, the matrix $I_n+A(p)$ is row stochastic and for all $k\in\N$ one has $p(k)\in(0,1)$.
\end{remark}

{\bf CODA model:}  provides the main model under study in this paper. This model assumes that each individual has access only to the action of its neighbors. The opinion of agent $i\in\V$ in this case updates according to the following rule:
\begin{align}
p_i(k+1)=p_i(k)+\frac{p_i(k)(1-p_i(k))}{n_i}\sum_{j\in N_i} (q_j(k)-p_i(k)). \label{eq:mn2}
\end{align}
\begin{remark}
For an agent $i\in\V$, both COCA and CODA models propose interaction weights depending only on the opinion $p_i$.
\end{remark}
We can emphasize a natural partition of $\V$ in two subsets $N^-(k)=\{i\in\V\mid q_i(k)=0\}$ and $N^+(k)=\{i\in\V\mid q_i(k)=1\}$. The main objective of this paper is to study how these sets evolve in time and what is the behavior of the opinions $p_i(k)$ inside these sets.\\
Throughout the paper we denote by $n^-(k)$ and $n^+(k)$ the cardinality of $N^-(k)$ and $N^+(k)$, respectively. Similarly, for an agent $i$ we denote by $N_i^-(k)=N_i\cap N^-(k)$ and $N_i^+(k)=N_i\cap N^+(k)$ and by $n_i^-(k)$ and $n_i^+(k)$ the cardinalities of these sets.

\subsection{Related model}\label{sec:martins}

\begin{figure}[!ht]
\centering
\includegraphics[clip=true,trim=7cm 11cm 7cm 11cm, scale=0.7]{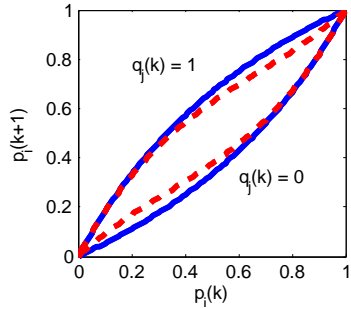}
\vspace{-0.3cm}
\caption{Comparison of the CODA model~\eqref{eq:mn2} and the model by Martins~\cite{Martins2008}. The two top curves display the update~\eqref{eq:mn2} and~\eqref{eq:MARTINS_CODA} when the influencing agent has action $q_j = 1$.
The two bottom curves display the update~\eqref{eq:mn2} and~\eqref{eq:MARTINS_CODA} when the influencing agent has action $q_j = 0$. Discontinuous red lines correspond to update~\eqref{eq:mn2} while strait blue lines correspond to update~\eqref{eq:MARTINS_CODA}.
}\label{fig:CODA_vs_Martins}
\end{figure}

The CODA model~\eqref{eq:mn2} was inspired by Martins' model~\cite{Martins2008} which was formulated in terms of the following bayesian update. Let us denote by $\tilde{p}_i(k)$ the opinion of agent $i\in\V$ at time $k\in\N$ when using Martins' model~\cite{Martins2008}. The updates of this opinion follows the rule described below. When agent $i$ is influenced by agent $j$,
\begin{equation}\label{eq:MARTINS_CODA}
\frac{\tilde{p}_i(k+1)}{1-\tilde{p}_i(k+1)} = 
\frac{\tilde{p}_i(k)}{1-\tilde{p}_i(k)} \cdot
\frac{\alpha}{1-\alpha}  \text{ if }  q_j(k) = 1,
\end{equation}
and where $\alpha$ is replaced by $1-\alpha$ if $q_j(k) = 0$. The constant $\alpha \ge 0.5$ is a model parameter linked to the amplitude of opinion change as a result of interactions. This parameter does not appear in our model. The study~\cite{Martins2008} is based on numerical experiments and do not present a theoretical analysis. The simulations found in~\cite{Martins2008} where obtained with $\alpha = 0.7$. The simulations appear to be qualitatively close to the ones resulting from our update~\eqref{eq:mn2} (see Section~\ref{sec:numerical_simulations}). To understand this fact, one can show that for $\alpha = 2/3 \approx 0.7$, and $q_j(k) = 1$, update~\eqref{eq:mn2} and~\eqref{eq:MARTINS_CODA} are equivalent for small $p_i(k)$ values. Similarly, if $q_j(k)=0$ the two updates are equivalent for $p_i(k)$ values close to $1$. For other $p_i(k)$ values, $p_i(k+1)$ still remains close to $\tilde{p}_i(k+1)$. These facts are illustrated in Figure~\ref{fig:CODA_vs_Martins}. As a consequence, our CODA model~\eqref{eq:mn2} can be seen as a consensus-type version of Martins' model. One advantage of our model is to explicitely explains why extremist agents with opinion close to $0$ or $1$ hardly change their actions. This is due to the weight $p_i (1-p_i)$ in update~\eqref{eq:mn2}.

\section{Analysis of the COCA model}\label{COCAmodel}

\sam{When agents have a direct access to opinions without quantification, opinions tend to a consensus at exponential speed as is shown in the present section. This is now longer the case when agents only have access the quantized actions instead of opinions of their neighbors, as will be shown in Section~\ref{CODAanalysis}.} It is noteworthy that an opinion stays constant according to dynamics \eqref{eq:mn1} if it is initialized at the value $0$ or $1$. Therefore, by the following assumption we exclude these extreme cases from the analysis.
\begin{assumption}\label{COCA1}
There exists a \textit{strictly} positive constants $ \epsilon\in\left(0,\frac{1}{2}\right)$ such that for all $ i\in\mathcal{V} $, one has $p_i^0\in[\epsilon,1-\epsilon]$.
\end{assumption}

For a subset of agents, $ \V'\subseteq \V $ and for all $ k\in \N $, let us define $m_{\V'}(k):=\min\limits_{j\in\V'}p_j(k) $ and $ M_{\V'}(k):=\max\limits_{j\in\V'}p_j(k)$. The following result can be easily proven by induction.
\begin{lemma} \label{COCA2}
Let $\mathcal{V}'\subseteq \mathcal{V}$ such that $\V'$ is isolated \ie for all $ i\in \mathcal{V}' $, for all $ j\in \mathcal{V}\setminus \mathcal{V}' $, $ (j,i)\notin \mathcal{E} $. Then,
\begin{enumerate}
\item $M_{\mathcal{V}'}(k+1)\leq M_{\mathcal{V}'}(k)$
\item $m_{\mathcal{V}'}(k+1)\geq m_{\mathcal{V}'}(k)$ 
\end{enumerate}
\end{lemma} 
Let us assume for the rest of this section that the graph is strongly connected. 

Applying Lemma \ref{COCA2} with $\V'=\V$ one obtains that sequences $m_{\V}(k)$ and $M_{\V}(k)$ are both monotonic and bounded, thus convergent. The interaction weight from the agent $j\in N_i$ to the agent $i$ at time $k$ is defined by $a_{ij}(k):=\frac{p_i(k)\big(1-p_i(k)\big)}{n_i}$. Taking into account that  $p_i^0\in[\epsilon,1-\epsilon]$, the behavior of the function $x\mapsto x(1-x)$ on $[0,1]$ and Lemma \ref{COCA2}, one deduces that $a_{ij}(k)\in[\frac{\epsilon(1-\epsilon)}{n_i},\frac{1}{4n_i}]$ for all $k\in\N$ and $j\in N_i$. Moreover, for all $i\in\V$ one has $1\leq n_i\leq n-1$ yielding \begin{equation}\label{aij}a_{ij}(k)\in\left[\frac{\epsilon(1-\epsilon)}{n-1},\frac{1}{4}\right], \quad \forall i\in\V,\ \forall j\in N_i.\end{equation}

Since $a_{ii}=1-\sum_{j\in N_i} a_{ij}$, straightforward computation shows that
 \begin{equation}\label{aii}a_{ii}(k)\in\left[\frac{3}{4},1-\epsilon(1-\epsilon)\right], \quad \forall i\in\V.\end{equation}
 
 Equation \eqref{aij} and \eqref{aii} shows that Assumption 1 in \cite{blondel2005} holds true. Moreover, in this paper we are dealing with a fixed (strongly) connected graph meaning that Assumption 2 and 3  in \cite{blondel2005} also hold. Therefore the following result is a consequence of Theorems 1 and Lemma 1 in \cite{blondel2005}.
 
 \begin{proposition}
 As far as the graph $\G$ is (strongly) connected and Assumption \ref{COCA1} is satisfied, the update rule \eqref{eq:mn1} guarantees asymptotic consensus of all opinions. Furthermore, the convergence rate to consensus satisfies
\[M_{\V}(k+1)-m_{\V}(k+1)\leq \beta(k) \Big(M_{\V}(k)-m_{\V}(k)\Big)\]
with $\beta(k)\leq\beta(\epsilon):=\Big(1-\min\{\frac{3}{4},\frac{ \epsilon(1-\epsilon)}{n-1}\}\Big)$.
 \end{proposition} 
 
 \section{Analysis of the CODA model}\label{CODAanalysis}
 
 For the remaining part of this work we consider the update rule \eqref{eq:mn2}. In this section we derive the set of possible equilibrium points and we give the conditions guaranteeing that the action $q_i$ of the agent $i\in\V$ is preserved/changed over time. Throughout the rest of the paper we assume that $\forall i\in\V, \ p_i(0)\in(0,1)$.
 
\subsection{Characterization of equilibria}

Let us define the following finite set of rational numbers:
\begin{equation}\label{S}
\S=\left\{\frac{k}{m} \ \Big|\ k,m\in\sam{\N}, k\leq m\leq n-1 \right\}\subset \Q .
\end{equation}
The main result of this section states that the possible equilibrium points of the opinions belong to $\S$.  Let us first introduce an instrumental result. For the rest of the paper we use the notation $r_i(k):=\displaystyle\frac{n_i^+(k)}{n_i}$.
\begin{lemma}\label{lemma:CODA}
Let $i\in\V$, $p_i(0)\in(0,1)$. Then, \sam{for all $k \in \N$}, one of the following relations holds
\begin{equation}\label{eq1}
p_i(k)< p_i(k+1)< r_i(k),\end{equation}
\begin{equation}\label{eq2}
p_i(k)>  p_i(k+1)> r_i(k),
\end{equation} 
\sam{or
\begin{equation}\label{eq_equal}
p_i(k) =  p_i(k+1) = r_i(k).
\end{equation} 
}
\sam{Moreover,} $\forall k\in\N,\ p_i(k)\in(0,1)$.
\end{lemma}  
\begin{IEEEproof}
The discrete dynamics \eqref{eq:mn2} can be rewritten by replacing $q_j(k)$ with $0$ or $1$ respectively. At time $k$ the agent $i$ has $n_i^-(k)$ neighbors having the action equals 0 and $n_i^+(k)$ with action equal 1. Then \eqref{eq:mn2}  rewrites as:
\begin{equation}\label{eq:model_rewritten_with_ni_plus}
\begin{array}{ll}
p_i(k+1)&=p_i(k)+\frac{p_i(k)(1-p_i(k))}{n_i}( n_i^+(k) -n_ip_i(k))\\
&=p_i(k)+p_i(k)(1-p_i(k))\left( r_i(k) -p_i(k)\right),
\end{array}
\end{equation}
where we have used the fact that $n_i^+(k)+n_i^-(k)$ is a constant equal to $n_i$.
\sam{We proceed by induction. Assume that $p_i(k) \in (0,1)$. Equation~\eqref{eq:model_rewritten_with_ni_plus} shows that} if $p_i(k)<\displaystyle r_i(k)$ then $p_i(k+1)> p_i(k)$ and reversely if $p_i(k)>\displaystyle r_i(k)$ then $p_i(k+1)< p_i(k)$. \sam{Finally, if $p_i(k) = \displaystyle r_i(k)$, $p_i(k+1) = p_i(k)$.}
Moreover,
\[
p_i(k+1)- r_i(k)=\left(p_i(k)- r_i(k)\right)(1-p_i(k)(1-p_i(k)))
\]
yielding $p_i(k+1)<\displaystyle r_i(k)$ if $p_i(k)<\displaystyle r_i(k)$\\ and  $p_i(k+1)>\displaystyle r_i(k)$ if $p_i(k)>\displaystyle r_i(k)$. In these cases, either \eqref{eq1} or \eqref{eq2} holds \sam{and $p_i(k+1) \in (0,1)$. Finally, if $p_i(k)=\displaystyle r_i(k)$, then equation~\eqref{eq_equal} holds true. With the assumption that $p_i(k) \in (0,1)$, we also obtain $p_i(k+1) \in (0,1)$.}
\end{IEEEproof}
\begin{proposition}\label{CODAeq}
\sam{Let $i\in\V$ and suppose Assumption \ref{COCA1} holds. If $\big(n_i^+(k)\big)_{k\geq0}$ (and thus  $\big(n_i^-(k)\big)_{k\geq0}$) is stationary, then the sequence of opinions $\big(p_i(k)\big)_{k\geq0}$ converges and has the following limit :
$$p_i^*=\displaystyle\lim_{k\rightarrow\infty}p_i(k)=\displaystyle\frac{\displaystyle\lim_{k\rightarrow\infty}n_i^+(k)}{n_i}\in\S.$$
Conversely, if the sequence of opinions $\big(p_i(k)\big)_{k\geq0}$ converges towards a value $p_i^* \in (0,1)$, then the sequence $\big(n_i^+(k)\big)_{k\geq0}$ is stationary and the equality between the limits still holds.}
\end{proposition}

\begin{IEEEproof}
If the sequence $\big(n_i^+(k)\big)_{k\geq0}$ is stationary, one gets that $n_i^+(k)$ is constant for $k$ bigger than or equal to a fixed integer $k^*$. Let us denote by $\sam{\rho}_i^*$ the value of $ r_i(k)$ for $k\geq k^*$. \sam{According to Assumption~\ref{COCA1}, $p_i(0) \in (0,1)$, so that Lemma~\ref{lemma:CODA} applies. By induction one has either 
\begin{equation}\label{eq:increase_of_pi}
\forall k\geq k^*, p_i(k)\leq p_i(k+1)\leq \sam{\rho}_i^*
\end{equation}
or  
\begin{equation}\label{eq:decrease_of_pi}
\forall k\geq k^*, p_i(k)\geq p_i(k+1)\geq \sam{\rho}_i^*.
\end{equation}}
In the first case $\big(p_i(k)\big)_{k\geq0}$ is increasing and upper-bounded, in the second $\big(p_i(k)\big)_{k\geq0}$ is decreasing and lower-bounded. Thus $\big(p_i(k)\big)_{k\geq0}$ converges. \sam{Denote $p_i^*$ its limit. If $p_i^* = 1$, equation~\eqref{eq:increase_of_pi} must be satisfied and $\rho_i^* = 1$. If $p_i^* = 0$, equation~\eqref{eq:decrease_of_pi} must be satisfied and $\rho_i^* = 0$. Finally, if $p_i^* \in (0,1)$, one can take the limit of equation~\eqref{eq:model_rewritten_with_ni_plus} to show that $p_i^* = \rho_i^*$.}

\sam{Reversely, assume that the sequence of opinions $\big(p_i(k)\big)_{k\geq0}$ converges to a limit $p_i^* \in (0,1)$. So that for some $k^*\in \N$, $\forall k\ge k^*, p_i(k) \in (0,1)$. As a consequence, equation~\eqref{eq:model_rewritten_with_ni_plus} rewrites as 
\begin{align*}
\frac{p_i(k+1) - p_i(k)}{p_i(k)(1-p_i(k))} + p_i(k) &=  r_i(k).
\end{align*}
Taking the limit of the previous equation shows that $(n_i^+(k)/n_i)_{k\ge 0}$ converges and is thus stationary and its limit satisfies $\rho_i^* = p_i^*$. Since $n_i^+(k), n_i  \in \N$ and $n_i^+(k) \le n_i$, it is clear that $p_i^*\in\S$.
}
\end{IEEEproof}

\subsection{Preservation of actions}

In this subsection we investigate the conditions ensuring that the action $q_i$ does not change over time. More precisely, we provide a criterion to detect when a group of agents sharing the same action will preserve it for all time. For the sake of simplicity we focus on $q_i(0)=0$ but using similar arguments the same results can be obtained for $q_i(0)=1$.

\begin{lemma}\label{lemma:actionpreserv}
Let $i\in\V$, if $n_i^-(k)\geq n_i^+(k)$ and $q_{i}(k)=0$ then $q_{i}(k+1)=0$.
\end{lemma}
\begin{IEEEproof}
Following the proof of Proposition \ref{CODAeq} one obtains that either \eqref{eq1}, \eqref{eq2} \sam{or \eqref{eq_equal} holds}. \\ If \eqref{eq1} \sam{or \eqref{eq_equal}} holds then $p_{i}(k+1)\leq\displaystyle r_i(k)$. Since  $n_i^-(k)\geq n_i^+(k)$  one deduces that $p_{i}(k+1)\leq\displaystyle\frac{1}{2}$ and since $q_{i}(k)=0$ one obtains that $q_{i}(k+1)=0$. \\
If \eqref{eq2} holds, then $p_{i}(k+1)\leq p_i(k)$. Since $q_{i}(k)=0$ one deduces that $p_i(k)\leq\displaystyle\frac{1}{2}$.  Therefore, $p_{i}(k+1)\leq\displaystyle\frac{1}{2}$ meaning again that $q_{i}(k+1)=0$.
\end{IEEEproof}

Throughout the paper we denote by $|A|$ the cardinality of a set $A$.

\begin{definition}\label{Polarized-Cluster}
We say a subset of agents $A\subseteq\V$ is a \textit{robust polarized cluster} if the following hold:
\begin{itemize}
\item $\forall i,j\in A,\  q_i(0)=q_{j}(0)$;\\

\item $\forall i\in A,\ |N_i\cap A|\geq |N_i\setminus A|$.
\end{itemize}
\end{definition}
\begin{remark}Notice that, in this section we do not make any assumption on the connectivity of the interaction graph. This means in particular that it may happen to have $n_i=0$ for some agents belonging to $A$ above. 
\end{remark}
The next result explains why the word \textit{robust} appears in the previous definition.

\begin{proposition}\label{prop:actionpreserv}
If $A$ is a robust polarized cluster with $q_i(0)=0$ for a certain $i\in A$ then 
\begin{itemize}
\item $\forall i\in A,\ \forall k\in\N, \ q_i(k)=0$;\\

\item $\forall i\in A, \ \displaystyle\lim_{k\rightarrow\infty}p_i(k)\leq\displaystyle\frac{1}{2}$.
\end{itemize}
\end{proposition}

\begin{IEEEproof}
The proof will be done by induction. Let us remark first that following the first item of Definition \ref{Polarized-Cluster} we have $q_i(0)=0,  \forall i\in A$.  Let us assume that for a fixed $k^*$ and $\forall i\in A$ one has $q_i(k^*)=0$. Let us also recall that the interaction graph under consideration is fixed. Therefore we still have that  $|N_i\cap A|\geq |N_i\setminus A|$. Moreover,  $|N_i\cap A|\subseteq N_i^-(k^*)$ implying that $n_i^-(k^*)\geq n_i^+(k^*), \forall i\in A$. Thus,  Lemma \ref{lemma:actionpreserv} yields that $ \forall i\in A$ one has $q_i(k^*+1)=0$ and the proof ends. 
\end{IEEEproof}

\sam{Proposition~\ref{prop:actionpreserv} is illustrated in Sections~\ref{sec:simu_minority} and~\ref{sec:simu_lattice}.}

\subsection{Change and diffusion of actions}

The goal of this subsection is twofold. First to provide conditions at a given time $k\in\N$ guaranteeing that the action of a fixed agent $i\in\V$ will change at time $k+1$. Secondly, we analyze the propagation/diffusion of the action of a robust polarized cluster inside the overall network. Due to the symmetry of reasonings we continue to focus only on one case $q_i(0)=0$ or $q_i(0)=1$.

\begin{proposition}\label{prop:changeopaction}
Let $i\in\V$ and $k\in\N$ such that $p_{i}(k)>\displaystyle\frac{1}{2}$ and $n_i^-(k)>n_i^+(k)$. Let $\epsilon(n_i^+(k),n_i) \in \left(0,\displaystyle\frac{1}{2}\right)$ be the unique positive real solution of the equation:
\begin{equation}\label{eps-x def}
x^3-\left(\frac{1}{2}- r_i(k)\right)x^2-\frac{3}{4}x+\left(\frac{1}{2}- r_i(k)\right)\frac{1}{4}=0.
\end{equation}  Then $p_i(k+1)\leq\displaystyle\frac{1}{2}$ if and only if $p_i(k)<\displaystyle\frac{1}{2}+\epsilon(n_i^+(k),n_i)$.
\end{proposition}
\begin{IEEEproof}
Let us introduce $x:=p_i(k)-\displaystyle\frac{1}{2}\in\left(0,\displaystyle\frac{1}{2}\right]$.
Recall that \eqref{eq:mn2}  rewrites as:
\begin{align*}
p_i(k+1)=p_i(k)+p_i(k)(1-p_i(k))\left(  r_i(k) -p_i(k)\right).
\end{align*}
Consequently, $p_i(k+1)\leq\displaystyle\frac{1}{2}$ if and only if 
\[
p_i(k)+p_i(k)(1-p_i(k))\left(  r_i(k) -p_i(k)\right)\leq\frac{1}{2},
\] which rewrites in term of $x$ as
\[
x^3+\left(\frac{1}{2}- r_i(k)\right)x^2+\frac{3}{4}x-\left(\frac{1}{2}- r_i(k)\right)\frac{1}{4}\leq0.
\] 
\sam{Since $n_i^+(k)/n_i < 1/2$,} the inequality above is \sam{strictly} satisfied for $x=0$ but it does not hold for $x=\displaystyle\frac{1}{2}$. Therefore, by continuity, \eqref{eps-x def} has a solution into $\left(0,\displaystyle\frac{1}{2}\right)$. Moreover, by using the first derivative it is straightforward that the function
\[x\mapsto x^3+\left(\frac{1}{2}- r_i(k)\right)x^2+\frac{3}{4}x-\frac{1}{4}\left(\frac{1}{2}- r_i(k)\right)\] is strictly increasing for $x\geq0$. Thus \eqref{eps-x def} has an unique positive solution Denoting it by $\epsilon(n_i^+(k))$ one deduces that $p_i(k+1)\leq\displaystyle\frac{1}{2}$ if and only if  $p_i(k)\in\left[0,\displaystyle\frac{1}{2}+\epsilon(n_i^+(k),n_i) \right)$ provided that $n_i^-(k)>n_i^+(k)$.
\end{IEEEproof}
The previous result states that an agent will change its action when it is influenced by more opposite actions, only if its opinion is sufficiently close to the boundary between the actions. The notion of sufficiently closed  depends on the proportion of opposite actions that influence the agent and is exactly quantified by $\epsilon(n_i^+(k),n_i)$.

\begin{lemma}\label{lemma:actionchange}
\sam{Let $i\in\V$ and $T_i\in\N$ such that $n_i^-(k)> n_i^+(k)$, $ \forall k\geq T_i$.} Then it exists $k^*\geq T_i$ such that $q_{i}(k)=0,\ \forall k\geq k^*$.
\end{lemma}

\begin{IEEEproof}
To obtain a contradiction, let us suppose that $p_i(k)\geq\displaystyle\frac{1}{2},\ \forall k\geq T_i$. Since $n_i^-(k)> n_i^+(k)$, $ \forall k\geq T_i $ we get that $\displaystyle r_i(k)<\displaystyle\frac{1}{2}, \forall k\geq T_i$.  We recall that
\[
p_i(k+1)=p_i(k)+p_i(k)(1-p_i(k))\left( r_i(k)-p_i(k)\right).
\]
Therefore, the sequence $\big(p_i(k)\big)_{k\geq T_i}$ is strictly decreasing and lower-bounded by $\displaystyle\frac{1}{2}$. Consequently $\big(p_i(k)\big)_{k\in\N}$ converges and its limit is bigger than or equal to $\displaystyle\frac{1}{2}$ \sam{and cannot be equal to $1$}. \\[2mm]
From Proposition \ref{CODAeq} the limit of $\big(p_i(k)\big)_{k\in\N}$ is $\displaystyle\frac{\displaystyle\lim_{k\rightarrow\infty}n_i^+(k)}{n_i}$.
On the other hand $n_i^-(k)> n_i^+(k)$, $ \forall k\geq T_i$ \sam{and these numbers are non-negative.}
Therefore $\displaystyle\frac{\displaystyle\lim_{k\rightarrow\infty}n_i^+(k)}{n_i}<\displaystyle\frac{1}{2}$ contradicting the fact that $\big(p_i(k)\big)_{k\in\N}$ converges toward a value bigger than or equal to $\displaystyle\frac{1}{2}$.\\[2mm]
The previous reasoning shows that  it exists $k^*\geq T_i$ such that $q_{i}(k^*)=0$. To get $q_{i}(k)=0,\ \forall k\geq k^*$ it is sufficient to recursively apply Lemma \ref{lemma:actionpreserv}.
\end{IEEEproof}

The next result characterizes the diffusion of the action of a robust polarized cluster over the network.

\begin{proposition} \label{prop:actionchange}
Let us consider the sets $A_1, A_2,\ldots, A_d$ such that  
\begin{itemize}
\item $A_1$ is a robust polarized set with $q_i(0)=0$ for a certain $i\in A_1$ (and thus $\forall i\in A_1,\ q_i(0)=0$);\\

\item $\forall h\in\{1,\ldots,d-1\}$ and $\forall i\in A_{h+1}$ one has 
$$|N_i\cap A_h|> |N_i\setminus A_h|.$$
\end{itemize}
\sam{Then, for all $h\in\{1,\ldots,d\}$, $\exists T_h\in\N$, such that, 
$$
\forall k \geq T_h, \forall i\in A_h, q_i(k)=0,
$$
where one can choose $T_1=0$.
}
Consequently,  $\forall i\in \displaystyle\bigcup_{h=1}^d A_h$ one has  $\displaystyle\lim_{k\rightarrow\infty}p_i(k)\leq\displaystyle\frac{1}{2}$. 
\end{proposition}
\begin{IEEEproof}
\sam{The proposition with $h=1$ and $T_1 = 0$ follows from the Proposition \ref{prop:actionpreserv}.}
\sam{For $h>1$,} we proceed recursively. 
Assume that the proposition holds for $h\in\{1,\ldots, f\}$ with $f<d$ and show that it holds for $h=f+1$. We know that $\forall i\in A_{f+1}$ one has 
$$|N_i\cap A_f|> |N_i\setminus A_f|.$$ Moreover $N_i\cap A_f\subseteq N_i^-(k),\ \forall k\geq T_f$ and $\forall i\in A_{f+1}$. Therefore we can apply Lemma \ref{lemma:actionchange} for any $ i\in A_{f+1}$. Choosing $T_{f+1}=\displaystyle\max_{i\in A_{f+1}}\ T_i$ one obtains that the proposition holds for $h=f+1$.

The last part of the statement is a simple consequence of the fact that $\forall i\in \displaystyle\bigcup_{h=1}^d A_h$ one has $q_i(k)=0, \forall k\geq T$ where $T=\displaystyle\max_{h\in \{2,\ldots,d\}}\ T_h$.
\end{IEEEproof}

\sam{Proposition~\ref{prop:actionpreserv} is illustrated in Section~\ref{sec:simu_lattice}.}

\section{Some particular network topologies}\label{Particular-Top}

\subsection{Complete graph}

In this subsection we use previous results to completely characterize the opinion dynamics when the interactions are described by the complete graph. 

\begin{proposition}\label{complete}
If $n^-(0)>n^+(0)$ then $\forall i\in\V$ the limit behavior of the opinion is given by $\displaystyle\lim_{k\rightarrow\infty}p_i(k)=0$. Reversely, $n^+(0)>n^-(0)$ implies $\displaystyle\lim_{k\rightarrow\infty}p_i(k)=1$.
\end{proposition}
\begin{IEEEproof}
If $n^-(0)>n^+(0)$, due to all-to-all connections one has that $N^-(0)$ is a robust polarized cluster. Moreover, $\forall i\in N^+(0)$ one has that $|N_i\cap N^-(0)|> |N_i\setminus N^-(0)|.$ Thus, applying Proposition \ref{prop:actionchange} with $A_1=N^-(0)$ and $A_2=N^+(0)$ we obtain that it exists a value $T\in\N$ such that $q_i(k)=0, \ \forall i\in\V, \forall k\geq T$. Following Proposition \ref{CODAeq} this yields that $$\displaystyle\lim_{k\rightarrow\infty}p_i(k)=0.$$
The case $n^+(0)>n^-(0)$ is proven in a similar way and as we have done throughout the paper we omit the details.
\end{IEEEproof}

Let us consider now the case $n^-(0)=n^+(0)$. It is noteworthy that in this case $n=|\V|$ is even and
\begin{itemize}
\item $\forall i\in N^+(0)$ one has $$|N_i\cap N^-(0)|=\frac{n}{2}>\frac{n}{2}-1= |N_i\setminus N^-(0)|.$$
\item $\forall i\in N^-(0)$ one has $$|N_i\cap N^+(0)|=\frac{n}{2}>\frac{n}{2}-1= |N_i\setminus N^+(0)|.$$
\end{itemize}
If the initial conditions are not symmetric w.r.t. $\displaystyle\frac{1}{2}$ an agent will cross from  $N^+$ to $N^-$ (or reversely) and we recover the situation treated in Proposition \ref{complete}. Therefore, to finish the analysis we give the following result that deals with initial conditions symmetrically displayed w.r.t. $\displaystyle\frac{1}{2}$. This case emphasizes an interesting oscillatory behavior of the opinions.

\begin{proposition}\label{completesym}
Assume that $n^+(0)=n^-(0)$ and moreover $\forall i\in\{1,\ldots,\frac{n}{2}\}$ there exist $\eta_i(0)\in\left(0,\displaystyle\frac{1}{2}\right)$ such that
\begin{equation}\label{init-cond complete-graph}
p_i(0)=\frac{1}{2}-\eta_i(0) \  \mbox{ and } p_{\frac{n}{2}+i}(0)=\frac{1}{2}+\eta_i(0).
\end{equation}
Then $n^+(k)=n^-(k),\ \forall k\in\N$ and $\exists k^*\in\N$ such that \sam{$\forall k\geq k^*, \forall j\in\V$,}
\[\begin{split}&|p_j(k)-\frac{1}{2}|\leq \epsilon^* \text{ and } \left(p_j(k)-\frac{1}{2}\right)\left(p_j(k+1)-\frac{1}{2}\right)<0,
\end{split}\]
 where $\epsilon^*<\frac{1}{6(n-1)}$ is the unique positive solution of the equation
\[
x^3+\frac{1}{2(n-1)}x^2+\frac{3}{4}x-\frac{1}{8(n-1)}=0.
\]
\end{proposition}
\begin{remark}
The result above states that all the agents in the network will finish by oscillating around $\frac{1}{2}$ in a $2\epsilon^*$ strip. 
\end{remark}
\begin{IEEEproof}
For all $i\in\{1,\ldots,\frac{n}{2}\}$ let us denote $\bar{i}=\frac{n}{2}+i$ and observe that due to all-to-all communications one has \[\frac{n_i^+(0)}{n_i}=\frac{1}{2}+\frac{1}{2(n-1)} \mbox{ and } \frac{n_{\bar{i}}^+(0)}{n_{\bar{i}}}=\frac{1}{2}-\frac{1}{2(n-1)}.\]
For all the couples $(i,\bar{i})$ the following holds
\begin{align*}
p_i(1)&=p_i(0)+p_i(0)(1-p_i(0))\left(\frac{n_i^+(0)}{n_i}-p_i(0)\right)\\
&=p_i(0)+(\frac{1}{4}-\eta_i(0)^2)\left(\frac{1}{2(n-1)}+\eta_i(0)\right)\\
p_{\bar{i}}(1)&=p_{\bar{i}}(0)+p_{\bar{i}}(0)(1-p_{\bar{i}}(0))\left(\frac{n_{\bar{i}}^+(0)}{n_{\bar{i}}}-p_{\bar{i}}(0)\right)\\
&=p_{\bar{i}}(0)-(\frac{1}{4}-\eta_i(0)^2)\left(\frac{1}{2(n-1)}+\eta_i(0)\right).
\end{align*} Therefore, $p_i(1)$ and $p_{\bar{i}}(1)$ remain symmetric w.r.t. $\displaystyle\frac{1}{2}$. Inductively one obtains that $p_i(k)$ and $p_{\bar{i}}(k)$ remain symmetric w.r.t. $\displaystyle\frac{1}{2}$. In other words if an agent $i\in\{1,\ldots,\frac{n}{2}\}$ changes its action, the agent $\bar{i}$ will also change its action and the changing will be in the opposite sense. Consequently, $n^+(k)=n^-(k),\ \forall k\in\N$. \\
Let us show now that $\exists k^*\in\N$ such that \[|p_j(k)-\frac{1}{2}|\leq \epsilon^*,\ \forall k\geq k^*, \forall j\in\V.
\]
Notice that for all $i\in\{1,\ldots,\frac{n}{2}\}$ one has $p_i(0)<\frac{1}{2}$, $n_i^+(0)>n_i^-(0)$, $p_{\bar{i}}(0)>\frac{1}{2}$ and $n_{\bar{i}}^-(0)>n_{\bar{i}}^+(0)$.  Thus, from the proof of Lemma \ref{lemma:actionchange} $\exists k_i^*\in\N$ such that $p_{\bar{i}}(k_i^*+1)<\frac{1}{2}$ and from the reasoning above $p_i(k_i^*+1)>\frac{1}{2}$. From Proposition \ref{prop:changeopaction} this change of actions happens if and only if $p_{\bar{i}}(k_i^*)\in\left(\frac{1}{2},\frac{1}{2}+\epsilon^*\right)$ with $\epsilon^*=\epsilon(\frac{n}{2}-1,n-1) $ the unique positive solution of 
\[
x^3+\frac{1}{2(n-1)}x^2+\frac{3}{4}x-\frac{1}{8(n-1)}=0.
\]
For $x=\frac{1}{6(n-1)}$ the expression above is positive and for $x=\frac{1}{8(n-1)}$ the expression is negative yielding that $\epsilon^* \in\left(\frac{1}{8(n-1)},\frac{1}{6(n-1)} \right)$. Let $x=\frac{1}{2}-p_i(k^*)\leq\epsilon^*$ meaning that
\begin{align*}
&x^3+\frac{1}{2(n-1)}x^2+\frac{3}{4}x-\frac{1}{8(n-1)}\leq0.
\end{align*}
One notes that $p_i(k^*+1)<\frac{1}{2}+\epsilon^*$ is equivalent to
\[
x^3+\frac{1}{2(n-1)}x^2+\frac{3}{4}x+\epsilon^*-\frac{1}{8(n-1)}>0,
\] which is always true since $x\geq0$ and $\epsilon^*-\frac{1}{8(n-1)}>0$. Therefore, once $p_i$ enters in the tube of radius $\epsilon^*$ around $\frac{1}{2}$ it never escapes and moreover Proposition \ref{prop:changeopaction} applies at each iteration implying
\[
\left(p_j(k)-\frac{1}{2}\right)\left(p_j(k+1)-\frac{1}{2}\right)<0, 
\ \forall k\geq k^*, \forall j\in\V.\]
\end{IEEEproof}

\sam{Proposition~\ref{prop:actionpreserv} is illustrated in Section~\ref{sec:simu_complete_graph}.}

\subsection{Ring graph}

Throughout this section we consider the particular configuration in which the interactions are described by an undirected graph  in which each vertex has exactly two neighbors. In the following we identify agent $n+1$ as agent $1$ and agent $0$ as agent $n$. For a precise representation of the graph we assume that $\forall i\in\{1,\ldots,n\}$ one has $N_i=\{i-1,i+1\}$. 
\begin{proposition}\label{prop:ring}
Under the ring graph topology the opinions dynamics \eqref{eq:mn2} leads to the following properties:
\begin{itemize}
\item the set $\S$ defined in \eqref{S} reduces to $\left\{0,\displaystyle\frac{1}{2},1\right\}$;
\item if $\exists i\in\V$ such that $q_i(0)=q_{i+1}(0)=0$ then $\{i,i+1\}$ is a robust polarized cluster (\ie $q_i(k)=q_{i+1}(k),\ \forall k\in\N$);
\item if $\forall i\in\V$ one has $q_i(0)=1-q_{i+1}(0)$ then 
\begin{itemize}
\item either the initial opinions are not symmetric w.r.t. $\displaystyle\frac{1}{2}$ and agents will change actions asynchronously leading to robust polarized sets $\{i-1,i,i+1\}$.
\item or $p_i(0)\in\left\{\frac{1}{2}-\sigma,\frac{1}{2}+\sigma\right\},\forall i\in\V$ and agents will change actions synchronously preserving $n^-(k)=n^+(k),\forall k\in\N$. Moreover, for $\sigma$ solving 
$$8\sigma^3+8\sigma^2+14\sigma-1=0$$ one has  $p_i(k)\in\left\{\frac{1}{2}-\sigma,\frac{1}{2}+\sigma\right\},\forall i\in\V,\forall k\in\N$
\end{itemize}\end{itemize}
\end{proposition}
\begin{IEEEproof}
The first item is a consequence of Proposition \ref{CODAeq}. The second item follows from Proposition \ref{prop:actionpreserv}. 
The third item follows the ideas in Proposition \ref{completesym}. Imposing $p_i(k)=\frac{1}{2}+\sigma$, $p_i(k+1)= \frac{1}{2}-\sigma$ and using the ring graph topology we obtain:
\begin{align*}
\frac{1}{2}-\sigma=\frac{1}{2}+\sigma-\left(\frac{1}{2}+\sigma\right)^2\left(\frac{1}{2}-\sigma\right)
\end{align*} or equivalently 
\begin{align*}
8\sigma^3+8\sigma^2+14\sigma-1=0.
\end{align*}
\end{IEEEproof}
\section{Numerical illustrations}\label{sec:numerical_simulations}

\subsection{Influencial minority}\label{sec:simu_minority}

\begin{figure}[!ht]
\begin{center}
\begin{tabular}{cc}
\textbf{(A)} &
\textbf{(B)} \\
 \includegraphics[clip=true,trim=7.5cm 6cm 7cm 5cm, scale=0.5]{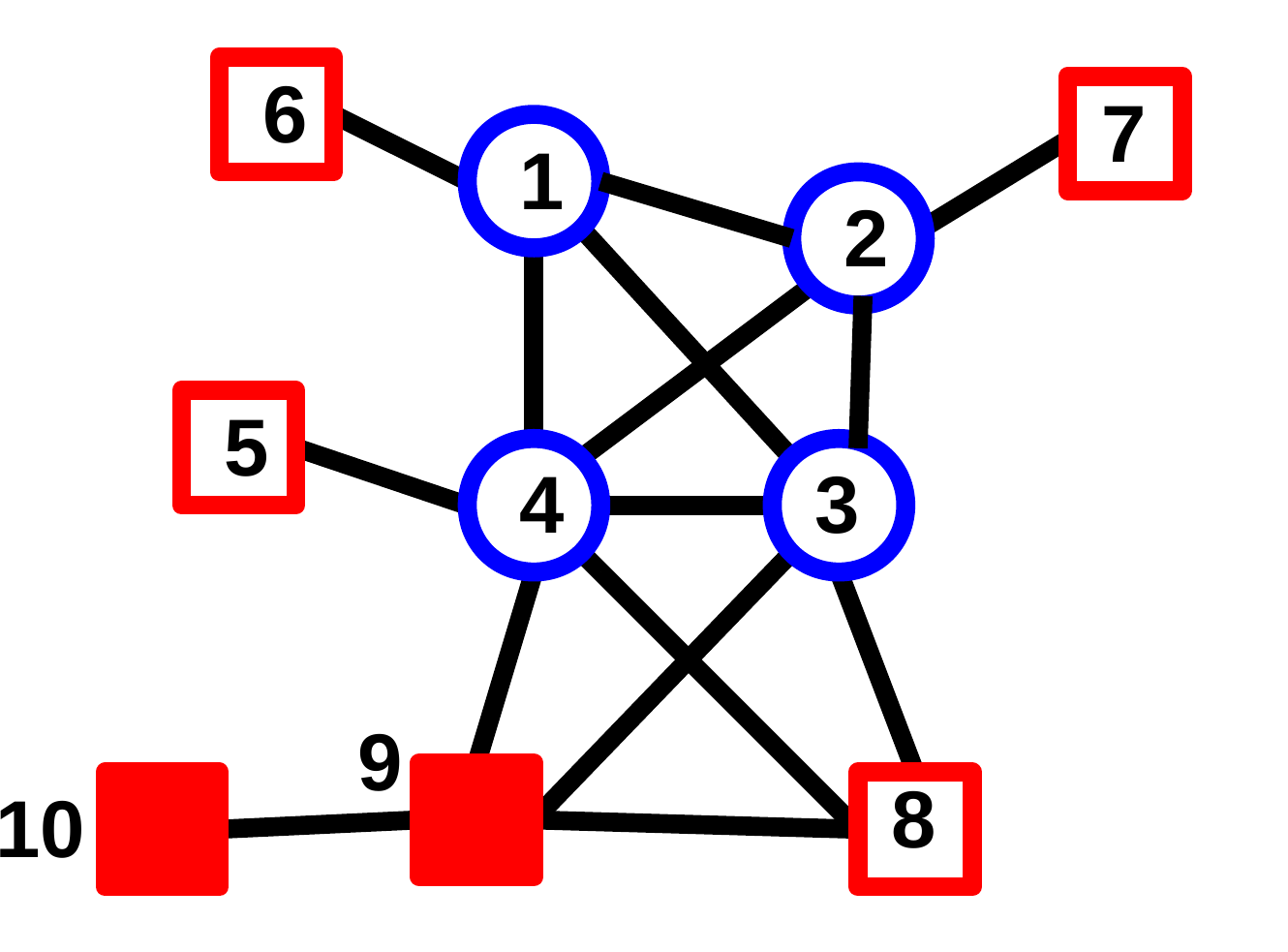} &
\includegraphics[clip=true,trim=3cm 8.5cm 3cm 9cm, scale=0.48]{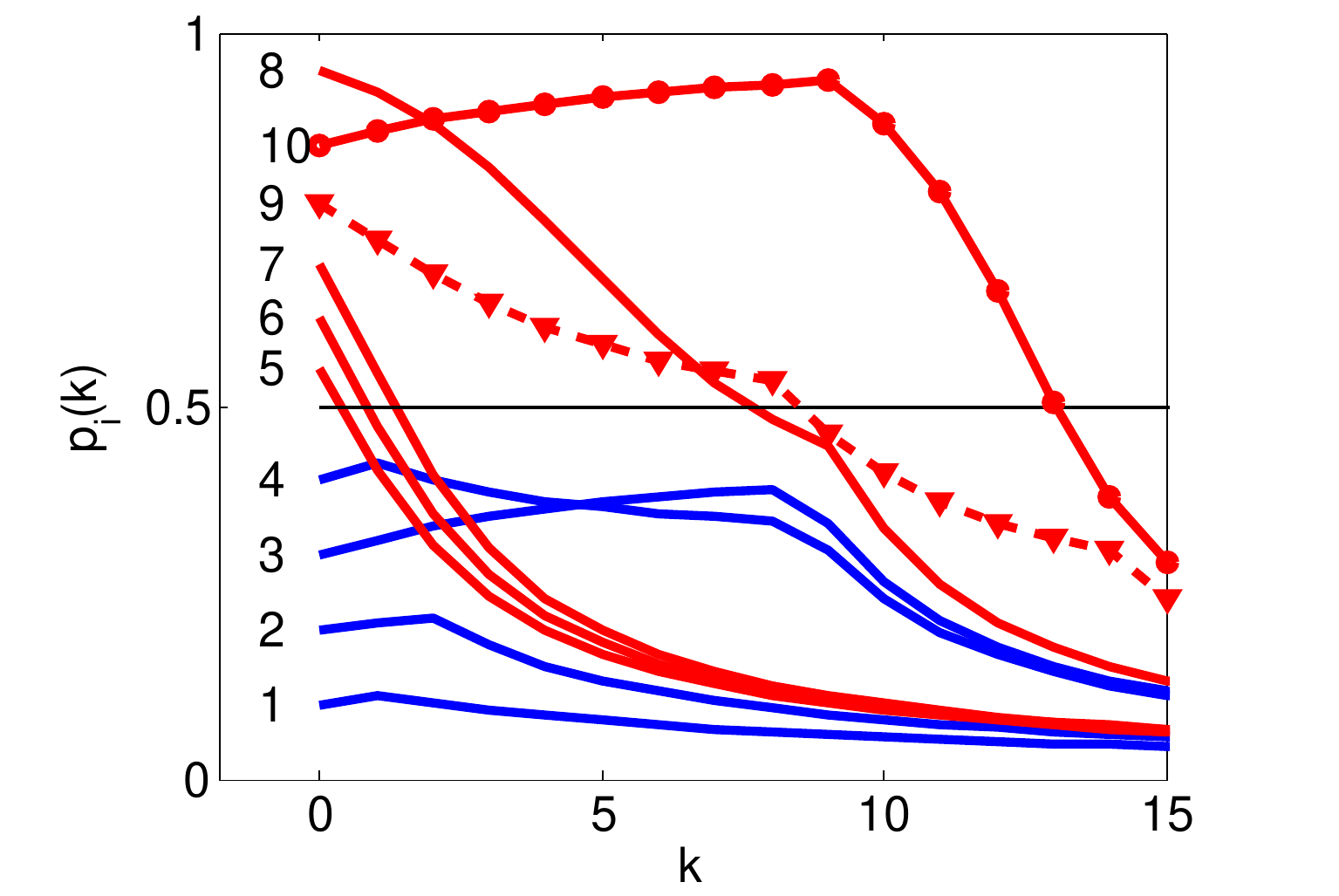}
\end{tabular}
\end{center}
\caption{Illustration of the action diffusion process described in Proposition~\ref{prop:actionchange}. Agents $1,2,3$ and $4$ start with $q_i(0) = 0$ and form a robust polarized cluster (Definition~\ref{Polarized-Cluster}), while agents $5,6,7,8,9$ and $10$ start with $q_i(0) = 1$. All agents converge to a state with action $\lim q_i = 0$.
}\label{fig:diffusion_of_action}
\end{figure}

We illustrate that a well connected minority can \textit{convince} a majority of agents located in the periphery of the interaction network. Proposition~\ref{prop:actionchange} can be used to predict the phenomenon given the topology of the social network and the initial actions only.
Figure~\ref{fig:diffusion_of_action} illustrates this fact. From Proposition~\ref{prop:actionchange}, taking $A_1 = \{1,2,3,4\}$, $A_2 = \{5,6,7,8\}$, $A_3 = \{9\}$ and $A_4 = \{10\}$, we predict that all agents will converge to a state with action $\lim q_i = 0$. We see in Figure~\ref{fig:diffusion_of_action}-B that initially agents $9$ and $3$ tend to approach $1/2$ since they have neighbors equally distributed over $1/2$. This behavior changes when agent $8$ passes the $1/2$ threshold.
Moreover, the action diffusion propagates to agent $10$ even though it originally had no neighbor with $q_i(0) = 0$. The decrease of agent $10$ towards $0$ only starts when agent $9$ passes the $1/2$ threshold.

\subsection{CODA on a square lattice}\label{sec:simu_lattice}

\begin{figure}[!ht]
\centering
\textbf{(A)} \hspace{2cm} \textbf{(B)} \hspace{1.8cm} \textbf{(C)} \hspace{0.5cm} \\
\includegraphics[clip=true,trim=5.4cm 8.5cm 6.8cm 9cm, scale=0.27]{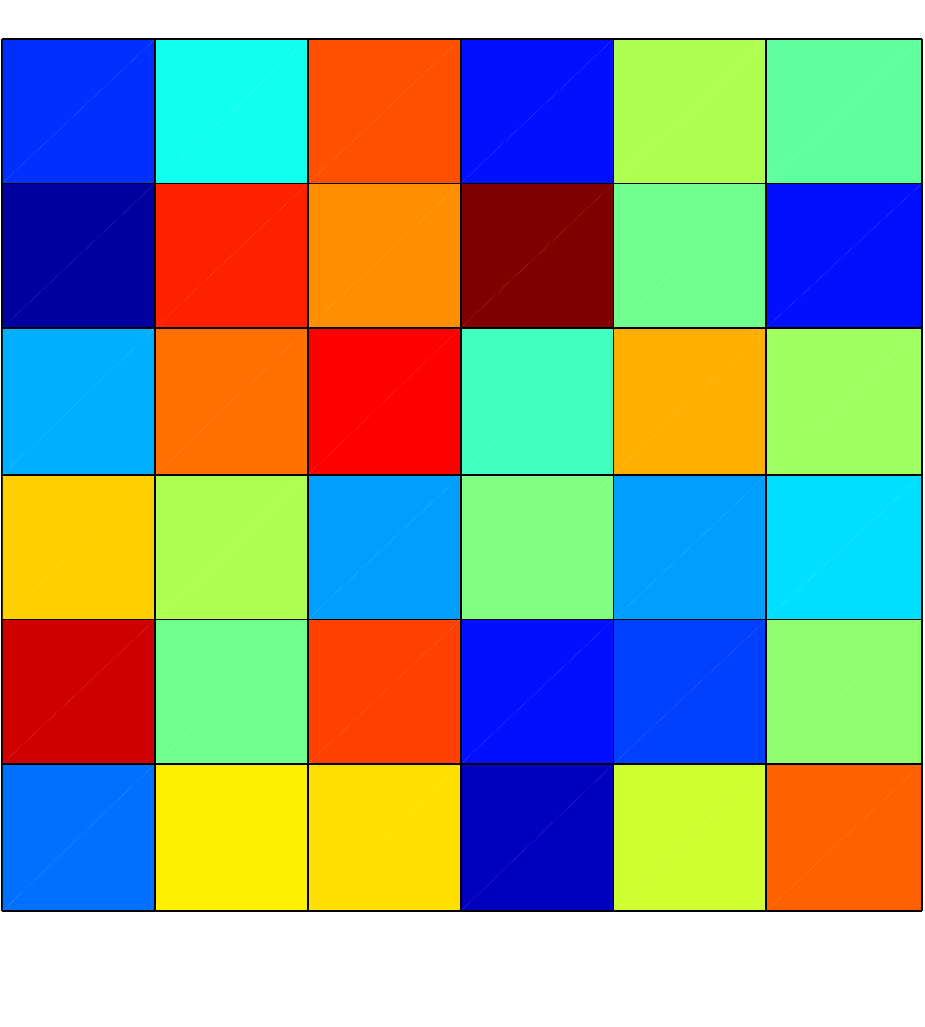}
\includegraphics[clip=true,trim=6cm 8.5cm 5.7cm 9cm, scale=0.27]{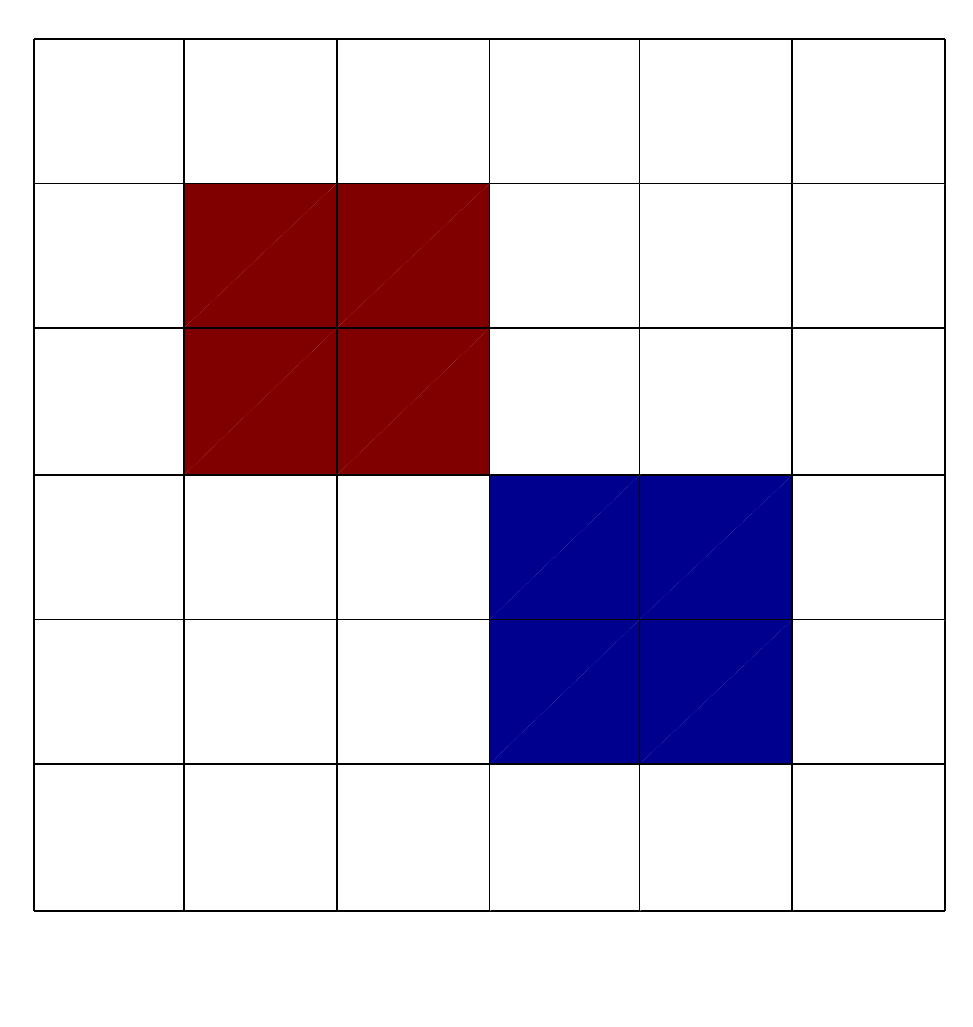}
\includegraphics[clip=true,trim=5.4cm 8.5cm 4.7cm 9cm, scale=0.27]{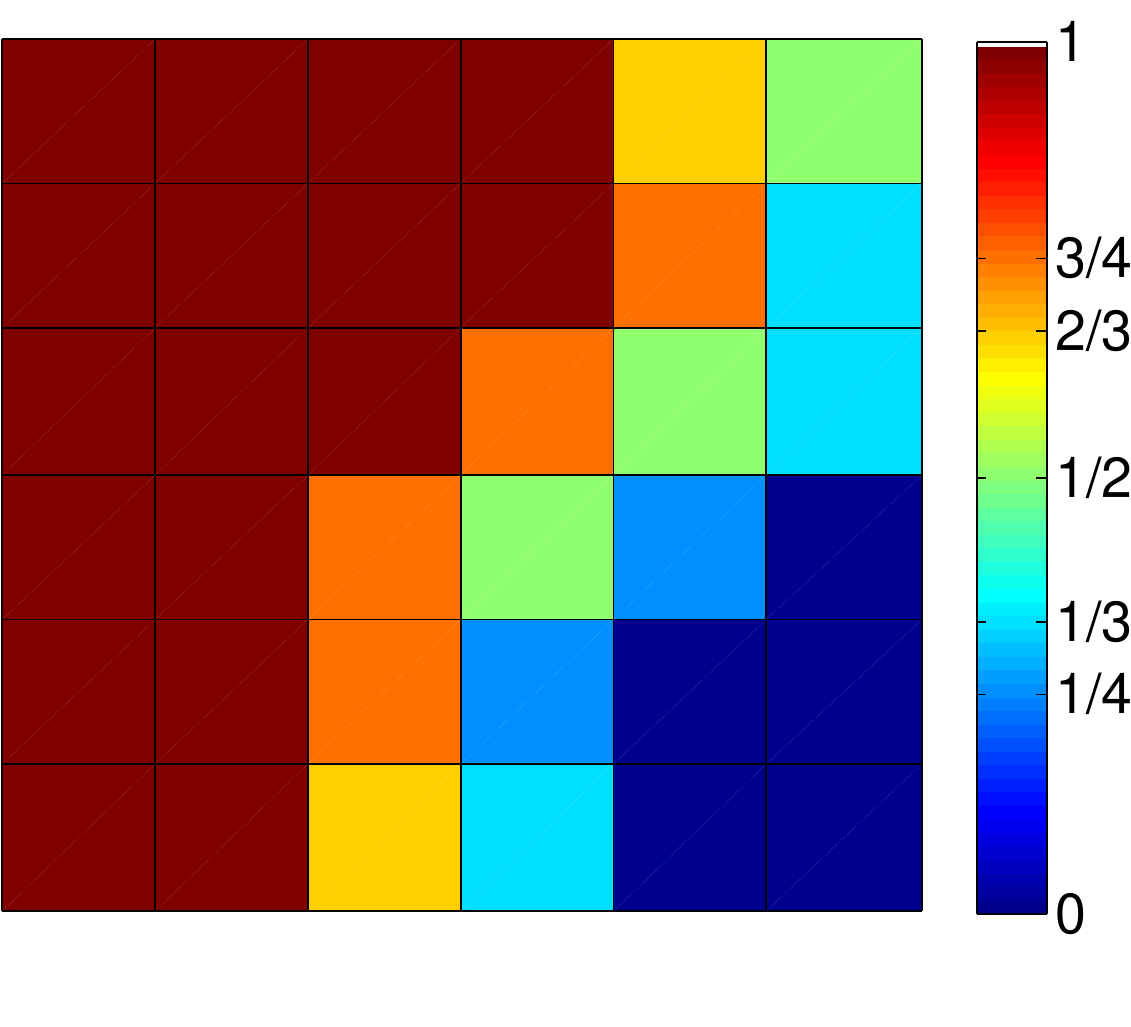}
\textbf{(D)} \\
\includegraphics[clip=true,trim=3cm 8.5cm 3cm 9cm, scale=0.5]{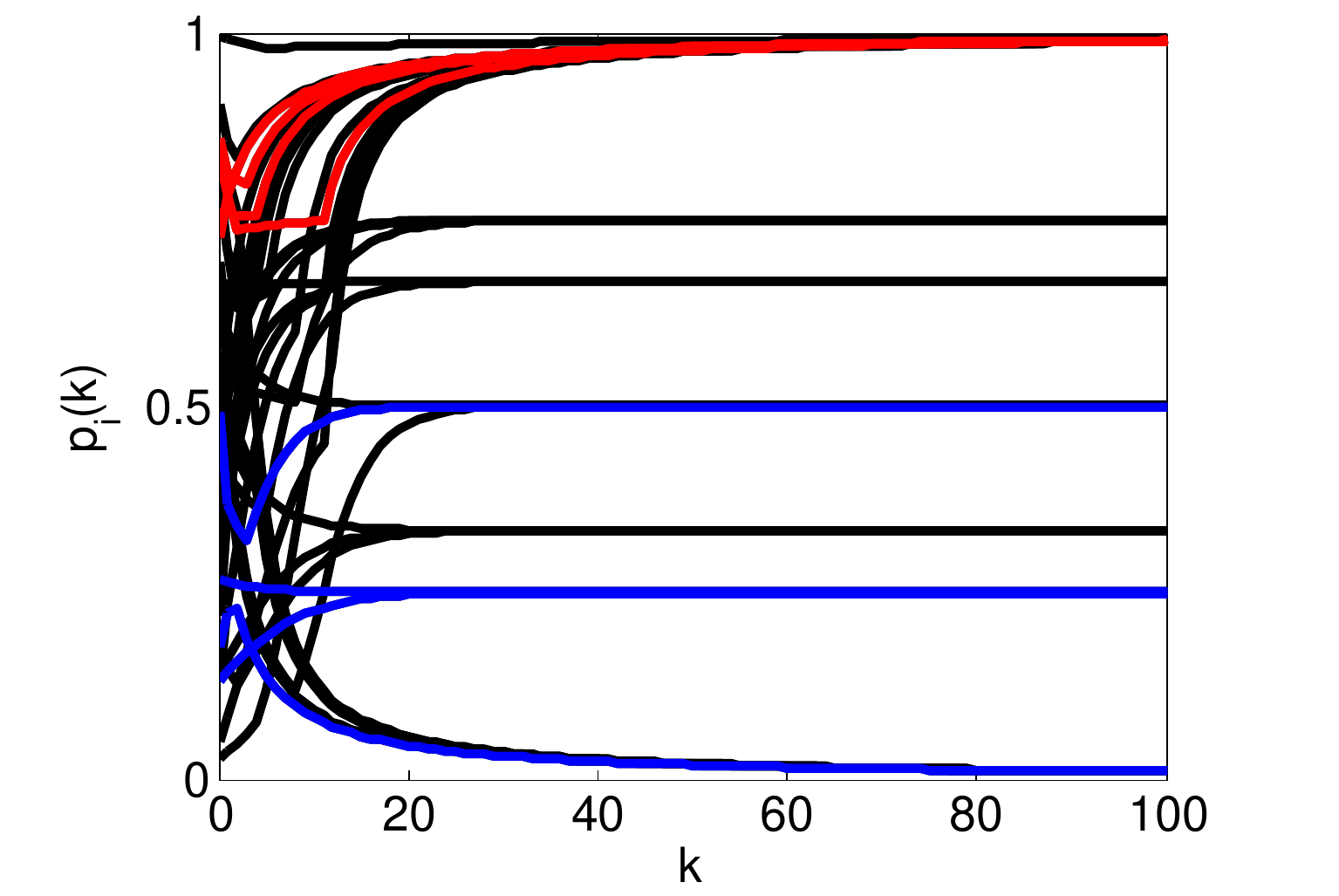}
\caption{Illustration of the CODA dynamics on a $6\times 6$ square lattice. (A) Each colored squared cell represents the initial opinion of an agent. Agents communicate with adjencent cells (above, below, left and right). (B) Robust clusters initially detected. (C) Final opinions of the agents. (D) Trajectories of the agents'opinions. Red lines correspond to agents belonging to a robust cluster with $q_i(0) = 1$. Blue lines correspond to agents belonging to a robust cluster with $q_i(0) = 0$.
}\label{fig:6-lattice}
\end{figure}

We illustrate our results when the topology of interactions is a square lattice. First, we use a $6\times 6$ lattice (see Figure~\ref{fig:6-lattice}). As illustrated in Figure~\ref{fig:6-lattice}-B, for this type of structure, the smallest robust clusters are formed by $2\times 2$ squares. As expected \sam{from Proposition~\ref{prop:actionpreserv}}, the robust clusters keep their inital actions and patches of same actions form around the robust clusters. Notice also that the values of convergence lie in set $\{0,1/4,1/3,1/2,2/3,3/4,1\}$ as predicted by Proposition~\ref{CODAeq}.

\begin{figure}[!ht]
\centering
\includegraphics[clip=true,trim=5.4cm 9.4cm 4.7cm 9cm, scale=0.7]{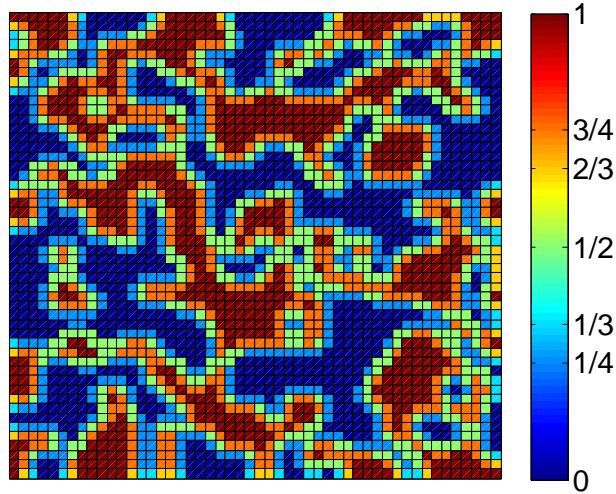}
\caption{Illustration of the CODA dynamics on a $50\times 50$ square lattice. Each colored squared cell represents the opinion of an agent after 100 iterations \sam{where initial conditions were drawn following independent uniform distributions}. Agents communicate with adjencent cells (above, below, left and right).
}\label{fig:50-lattice}
\end{figure}

Patches of agents with same actions are also observed for bigger lattice (see Figure~\ref{fig:50-lattice} for an instance final actions in a $50\times 50$ lattice \sam{where initial conditions were drawn following independent uniform distributions.}). This is in accordance with the patterns found in~\cite{Martins2008}.

\subsection{Oscillatory dynamics on a ring graph}\label{sec:simu_complete_graph}

\begin{figure}[!ht]
\centering
\includegraphics[scale=0.4]{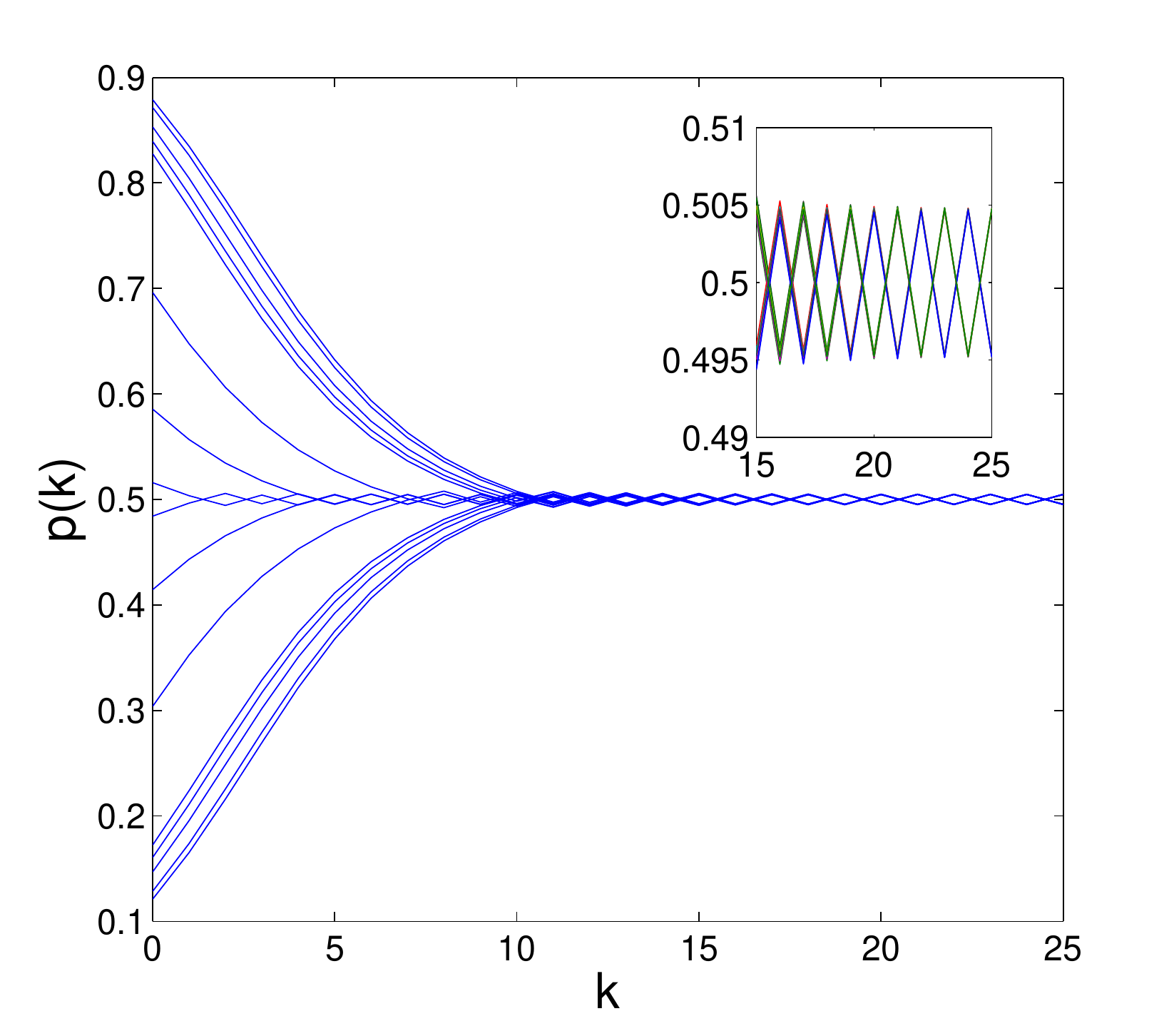}
\caption{Trajectories of $ 100 $ agents with the complete interaction graph and initial opinions distributed symmetrically around $1/2$.}
\label{fig:100-complete-graph}
\end{figure}

The following simulation displays the oscillation of agents'opinion around $1/2$ when the interaction graph is complete and when the initial opinions are symmetrically distributed around $1/2$ (see Proposition~\ref{completesym}). Figure~\ref{fig:100-complete-graph} shows an instance of this phenomenon for a system of $100$ agents.

\section{Conclusions}

In this paper we have introduced a novel opinion dynamics model in which agents have access to actions which are quantized version of the opinions of their neighbors. The model reproduces different behaviors observed in social networks such as dissensus, clustering, oscillations, opinion propagation. The main results of the paper provides the characterization of preservation and diffusion of action under general communication topologies. A complete analysis of the opinions behavior is given in the particular cases of complete and ring communication graphs. Numerical examples illustrate the main features of this model.

\bibliographystyle{IEEEtran}
\bibliography{diss}

\end{document}